%% file: lamhyparxiv.tex
\documentclass[12pt,a4paper,leqno]{amsart} 
\ifx\pdftexversion\undefined
  \usepackage[dvips]{graphicx,color}
\else
  \usepackage[pdftex]{graphicx,color}
\fi
\usepackage[francais]{babel}
\usepackage[latin1]{inputenc}%

\let\ifanglais\iffalse

\def\R{{\mathbb R}}
\def\pP{{\mathbb P}}

\def\vole{{\rm vol}_e}
\newtheoremstyle{mesthm}
  {10pt plus 1pt minus 1pt}
  {9pt minus 6pt}
{\slshape}
  {0.5cm}
  {\bfseries}
  {.}
  {1ex}
  {}
\newtheoremstyle{mesdefi}
  {6pt plus 1pt minus 1pt}
  {6pt plus 1pt minus 1pt}
  {}
  {0.5cm}
  {\bfseries}
  {.}
  {1ex}
  {}%

\theoremstyle{mesthm}
\newtheorem{lema}{\ifanglais{\large L}emma\else{\large L}emme\fi}
\newtheorem{theo}[lema]{\ifanglais{\large T}heorem\else {\large
    T}héorème\fi}

\newtheorem{prop}[lema]{{\large P}roposition} 
 
\newtheorem{cor}[lema]{\ifanglais{\large C}orollary\else{\large C}orollaire\fi}
\newtheorem{ppt}[lema]{\ifanglais{\large P}roperty\else{\large
    P}ropriété\fi}

\newtheorem{rmq}[lema]{\ifanglais{\large R}emark\else{\large
    R}emarque\fi}
\newtheorem{lemme}[lema]{\ifanglais{\large L}emma\else{\large L}emme\fi}

\theoremstyle{mesdefi}
\newtheorem{defi}[lema]{\ifanglais{\large D}efinition\else{\large
    D}éfinition\fi} 
 \newtheorem{rem}[lema]{\ifanglais{\large
    R}emark\else{\large  R}emarque\fi}
\newtheorem*{nb}{Remarque}

\makeatletter
\renewenvironment{proof}[1][{\sl \ifanglais Proof\else Démonstration\fi}]{\par
  \normalfont
  \topsep6\p@\@plus6\p@ \trivlist
  \item[\hskip\labelsep\slshape
    #1\@addpunct{.}]\ignorespaces
}{%
  \qed\endtrivlist
}
\makeatother

\newcommand{\n}[1]{\Vert #1 \Vert}

\title[$\lambda_1$ et $\delta$-hyperbolicité en géométrie de Hilbert]
{Bas du spectre et delta-hyperbolicit\'e en g\'eom\'etrie de 
Hilbert plane.}

\author[B.~Colbois et C.~Vernicos]{Bruno Colbois et Constantin Vernicos}
\address{Institut de math\'ematique\\
  Université de Neuch\^atel\\
  Rue \'Emile Argand 11\\
  CH--2007 Neuch\^atel\\
  Switzerland}
\email{Bruno.Colbois@unine.ch\\
Constantin.vernicos@unine.ch}

\subjclass[2000]{Differential Geometry; Metric Geometry}

\begin{document}

\maketitle
\begin{abstract}
  On montre l'\'equivalence entre l'hyperbolicit\'e au sens de Gromov 
  de la g\'eom\'etrie
de Hilbert d'un domaine convexe du plan et la non nullit\'e du bas du spectre de ce domaine.
\end{abstract}
\renewcommand{\abstractname}{Abstract}
\begin{abstract}(\textsl{Bottom of the spectrum and delta hyperbolicity in
Hilbert Plane geometry})
We prove that the Hilbert geometry of a convex domain in the plane is Gromov hyperbolic,
if, and only if, the bottom of its spectrum is not zero.
\end{abstract}

\input{spectrearxiv}

\nocite{*}

\bibliographystyle{amsalpha}
\bibliography{airhyp}

\end{document}

%% file: spectrearxiv.tex
\section*{Introduction}

Le but de cet article est de montrer l'\'equivalence entre l'hyperbolicit\'e 
au sens de Gromov de la g\'eom\'etrie
de Hilbert d'un domaine convexe $\mathcal C$ du plan et la non nullit\'e du bas du 
spectre $\lambda_1(\mathcal{C})$ de ce domaine.

En g\'eom\'etrie riemannienne, il existe des relations tr\`es fortes 
entre le bas du spectre du laplacien d'une vari\'et\'e compl\`ete de volume 
infini et la g\'eom\'etrie de cette vari\'et\'e. Par exemple, on sait 
que le bas du spectre d'une vari\'et\'e de Cartan-Hadamard \`a 
courbure sectionnelle $K \leq C<0$ est strictement positif. 
R\'ecemment, J. Cao \cite{jcao} a \'etudi\'e le cas des vari\'et\'es riemanniennes 
hyperboliques au sens de Gromov poss\'edant un quasi-p\^ole, et 
montr\'e que leur constante de Cheeger (donc le bas de leur spectre) 
\'etait strictement positif. Dans cet article, nous abordons ce type 
de questions dans le contexte des g\'eom\'etries de Hilbert.

Avant d'\'enoncer les r\'esultats pr\'ecis, rappelons qu'une g\'eom\'etrie de Hilbert
$(\mathcal C, d_{\mathcal C})$ est la donn\'ee d'un ouvert convexe et born\'e 
$\mathcal C$ de $\mathbb R^n$ muni de la distance de Hilbert $d_{\mathcal C}$ d\'efinie
de la mani\`ere suivante : 
pour toute paire de points distincts $p$ et $q$ dans $\mathcal{C}$, la droite passant
par $p$ et $q$ rencontre le bord $\partial \mathcal C$ de 
$\mathcal C$ en deux points distincts $a$ et $b$ tels que la droite passe
par $a$, $p$, $q$ et $b$ dans cet ordre. 

\begin{figure}[htpb]
\includegraphics{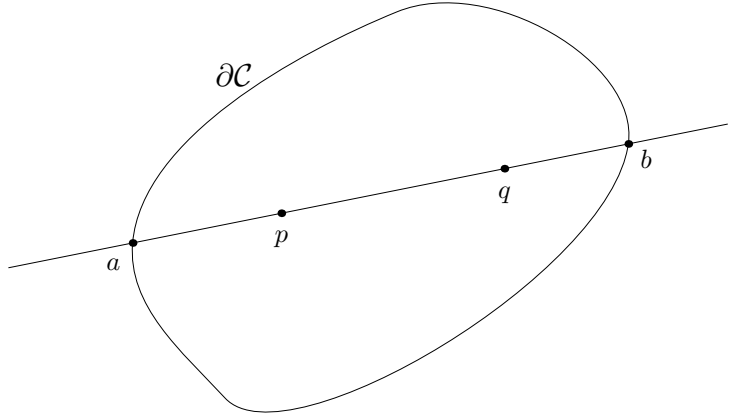}
\caption{Distance de Hilbert}
\end{figure}

On d\'efinit alors
$$d_{\mathcal C}(p,q) = \frac{1}{2} \ln [a,p,q,b]$$
o\`u $[a,p,q,b]$ est le birapport de $(a,p,q,b)$:
$$[a,p,q,b] = \frac{\Vert q-a \Vert_{e} }{\Vert p-a \Vert_{e}} \times 
\frac{\Vert 
p-b \Vert_{e}}{\Vert q-b \Vert_{e}}   $$
o\`u $ \Vert . \Vert_{e}$ d\'esigne la norme euclidienne.
On pose \'egalement $d_{\mathcal C}(p,p) =0$.

Sur $\mathcal C$, on peut mettre une norme de Finsler $C^0$, 
not\'ee  $\n{\cdot}_{\mathcal C}$, en proc\'edant
comme suit: si $p \in \mathcal C$ et $u_{p} \in T_{p}\mathcal C = \mathbb 
R^n$, la droite passant par $p$ et dirig\'ee par $u_{p}$ coupe $\partial \mathcal 
C$ en deux points $p^{+}$ et $p^-$.  
On pose alors
$$\n{u_{p}}_{\mathcal C} = \frac{1}{2} \n{u_{p}}_{e} 
\biggl(\frac{1}{\Vert 
p-p^{+} \Vert_{e}} + 
\frac{1}{\Vert p-p^- \Vert_{e}}\biggr) {\rm }$$

o\`u $ \n{u_{p}}_{e}$ d\'esigne la norme euclidienne de $u_{p}$.

Notons que la distance de longueur induite sur $\mathcal C$ par la norme 
$\Vert . \Vert_{\mathcal C}$ co\"\i ncide avec $d_{\mathcal C}$, mais 
nous n'utiliserons pas ce fait dans la suite.

A cette norme de Finsler est associ\'ee une norme duale: si $l_{p}$ 
est une forme lin\'eaire sur $T_{p}\mathcal C$, on pose
$$\n{l_{p}}_{\mathcal C}^{*} = \sup \{l_{p}(u_{p}): 
u_{p}\in T_{p}\mathcal C,\ \n{u_{p}}_{\mathcal C}=1 \}.$$

Gr\^ace \`a la norme de Finsler, on construit une forme volume et 
une mesure sur 
$\mathcal C$ (qui correspond en fait \`a  
la mesure de Hausdorff, voir \cite{bbi} 
exemple 5.5.13).

Soient 
$p \in \mathcal C$ et 
$$TB_{\mathcal C}(p) = \{u_{p} \in T_{p}\mathcal 
C =\mathbb R^n : 
\n{u_{p}}_{\mathcal C} < 1 \}$$
la boule unit\'e de $T_{p}\mathcal C$. 
Soit $\omega_{n}$ le volume de la boule unit\'e de l'espace euclidien
$\mathbb R^n$ 
et consid\'erons la fonction $h\colon\mathcal C \longrightarrow \mathbb R$ donn\'ee par 
${h(p)=\omega_{n}/\vole\bigl(B_{\mathcal C}(p)\bigr)}$ 
o\`u $\vole$ est le volume euclidien usuel.

Alors la mesure $\mu_{\mathcal C}$ (que nous nommerons   
mesure de Hilbert de $\mathcal{C}$, elle est \'egalement connue sous le nom de mesure
de Busemann)  associ\'ee \`a $\n{\cdot}_{\mathcal C}$ 
est d\'efinie ainsi : si $A \subset \mathcal C$ est un bor\'elien, on pose
$$\mu_{\mathcal C} (A) = \int_{A} h(p)d\vole(p)$$
o\`u $d\vole(p)$ est la mesure de Lebesgue.

Enfin, si $\Omega$ est un domaine avec $\overline{\Omega} \subset\mathcal C$, l'int\'egrale 
sur $\Omega$ par 
rapport \`a $\mu_{\mathcal C}$ d'une 
fonction $f$ d\'efinie sur $\Omega$ sera not\'ee
$\int_\Omega f~d\mu_\mathcal{C}.$

\begin{nb}
  Les r\'esultats de cet article restent inchang\'es si on consid\`ere une mesure
\'equivalente. En particulier, la mesure dite de Holmes-Thompson, qui
est \'equivalente par les in\'egalit\'es de Santalo et Bourgain-Milman \`a la mesure de Hilbert, donne les
les m\^emes r\'esultats.
\end{nb}

Lorsque $\mathcal{C}$ est une ellipse, $(\mathcal C, d_{\mathcal C})$ 
est le mod\`ele projectif (ou mod\`ele de Klein) de
la g\'eom\'etrie hyperbolique, et on peut penser aux g\'eom\'etries de Hilbert 
$(\mathcal C, d_{\mathcal C})$ comme \`a une g\'en\'eralisation naturelle 
de l'espace hyperbolique.
Une question communes \`a de nombreux travaux r\'ecents
(voir \cite{so},\cite{so2}, \cite{benoist},
\cite{cv}, \cite{cvv}, \cite{kn} et leurs r\'ef\'erences) est de d\'eterminer les propri\'et\'es
de l'espace hyperbolique dont ces g\'eom\'etries h\'eritent et de trouver
des caract\'erisations de l'espace hyperbolique parmi les g\'eom\'etries
de Hilbert. En particulier, dans \cite{benoist}, Y. Benoist obtient en 
toute dimension une caract\'erisation des convexes dont la 
g\'eom\'etrie de Hilbert associ\'ee est hyperbolique au sens de 
Gromov (voir \cite{brihae} pour une discussion de ce concept). Cette 
caract\'erisation est donn\'ee en fonction de la r\'egularit\'e du 
bord des convexes consid\'er\'es.

Dans cet article, on va caract\'eriser en dimension deux l'hyperbolicit\'e au sens de 
Gromov d'un point de vue spectral:
on va 
montrer qu'en dimension deux avoir un bas du spectre strictement positif est 
\'equivalent \`a \^etre hyperbolique au sens de Gromov.

On d\'efinit le bas du spectre de $\mathcal C$, que l'on note 
$\lambda_1(\mathcal C)$, par analogie avec ce qui se fait dans le cas 
des vari\'et\'es riemanniennes de volume infini. On pose
\begin{equation}
\label{eqdef1}
\lambda_1(\mathcal{C})=
\inf \frac{\displaystyle{\int_\mathcal{C}{\n{df_p}_\mathcal{C}^*}^2~d\mu_\mathcal{C}(p)}}
{\displaystyle{\int_\mathcal{C} f^2(p)d\mu_\mathcal{C}(p)}}
\end{equation}
o\`u l'infimum est pris sur toutes les fonctions 
lipschitziennes, non nulles, \`a support compact dans $\mathcal{C}$ et o\`u
$\mu_\mathcal{C}$ est la mesure de Hilbert associ\'ee \`a $\mathcal C$. 
 L'expression ci-dessus est appel\'ee le quotient de Rayleigh de $f$.

Le bas du spectre $\lambda_1(\mathcal{C})$ est un nombre r\'eel 
positif ou nul. On sait que lorsque $\mathcal C$ est une ellipse, 
c'est-\`a-dire que l'on se trouve dans le mod\`ele hyperbolique, il 
vaut $1/4$. On verra un peu plus loin que dans le cas o\`u 
$\mathcal C$ est un triangle, le bas du spectre est nul.

Notre r\'esultat principal est:
\begin{theo}\label{equiv1}
Soit $(\mathcal{C},d_\mathcal{C})$  un convexe du plan muni de sa m\'etrique de Hilbert.
Alors son bas du spectre est non nul, i.e., $\lambda_1(\mathcal{C})\neq 0$,
si et seulement s'il existe $\delta>0$ tel que $(\mathcal{C},d_\mathcal{C})$ est $\delta$-hyperbolique.
\end{theo}

Cela se g\'en\'eralise en dimension $n$ comme suit :
\begin{theo}
Soit $(\mathcal{C},d_\mathcal{C})$  un convexe de $\R^n$ muni de sa m\'etrique de Hilbert.
Alors il existe un r\'eel strictement positif $\lambda$ minorant 
le bas du spectre 
de toute section planaire, i.e. pour tout $H$ plan affine de $\R^n$ tel que $\mathcal{C}\cap H\neq \emptyset$
on a $$ \lambda_1(\mathcal{C}\cap H)\geq \lambda \text{,} $$
si et seulement s'il existe $\delta>0$ tel que $(\mathcal{C},d_\mathcal{C})$ est $\delta$-hyperbolique.
\end{theo}

Pour un domaine ouvert, born\'e (relativement \`a $d_{\mathcal C}$) $\Omega$ \`a bord lipschitzien de $\mathcal{C}$, 
on appellera bas
du spectre de Dirichlet le nombre 
\begin{equation}
\label{eqdef2}
\lambda_1^\mathcal{C}(\Omega)=
\inf \frac{\displaystyle{\int_\Omega{\n{df_p}_\mathcal{C}^*}^2~d\mu_\mathcal{C}(p)}}%
{\displaystyle{\int_\Omega f^2(p)d\mu_\mathcal{C}(p)}}\text{,}
\end{equation}
o\`u l'infimum est pris sur toutes les fonctions lipschitziennes, non nulles, 
\`a support dans $\Omega$.

Il est clair que $\lambda_1^\mathcal{C}(\Omega)\geq\lambda_1(\mathcal{C})$ et que $\lambda_1^\mathcal{C}$ est une fonction
d\'ecroissante suivant l'inclusion des ensembles: en effet, si 
$\Omega_{1} \subset \Omega_{2}$, toute fonction \`a support dans 
$\Omega_{1}$ s'\'etend naturellement par $0$ en une fonction \`a 
support dans $\Omega_{2}$.

L'id\'ee de la preuve du th\'eor\`eme \ref{equiv1} est la suivante. \`A la section 1, on montrera 
que la non nullit\'e du bas du spectre entra\^\i ne l'hyperbolicit\'e au 
sens de Gromov en nous appuyant principalement sur un r\'esultat de 
Y.~Benoist qui donne une condition suffisante pour qu'une famille de 
convexes munis de leur m\'etrique de Hilbert soit 
$\delta$-hyperbolique (proposition \ref{propbenoist} ci-dessous). On 
montrera que tous les convexes dont le bas du spectre est sup\'erieur 
ou \'egal \`a une constante positive $\lambda$ donn\'ee sont 
$\delta$-hyperboliques, $\delta$ d\'ependant de $\lambda$. 

\`A la section 2, on montrera la r\'eciproque en utilisant une 
in\'egalit\'e isop\'erim\'etrique que v\'erifie les espaces 
$\delta$-hyperboliques et sur une estimation du volume des boules, 
int\'eressante pour elle-m\^eme, valable en toute dimension, donn\'ee \`a la section \ref{volboule}. En 
particulier, le volume d'une boule de rayon fix\'e est uniform\'ement major\'e sur 
toutes les g\'eom\'etries de Hilbert de dimension $n$. Dans le cas 
riemannien cela est impliqu\'e par la donn\'ee d'un minorant sur la 
courbure de Ricci, mais on sait qu'une hypoth\`ese de ce type n'est 
pas v\'erifi\'ee par les g\'eom\'etries de Hilbert qui ne sont pas, en 
g\'en\'eral, des espaces d'Alexandroff. Pour conclure, nous aurons 
besoin d'une in\'egalit\'e du type Cheeger reliant les constantes 
isop\'erim\'etriques et le bas du spectre. Pour l'essentiel, on peut 
adapter ce qui se fait dans la cas riemannien en faisant un d\'etour 
par le cas finsl\'erien, mais cela pose quelques probl\`emes techniques. 
D'une part, les convexes que l'on consid\`ere poss\`edent une 
m\'etrique de Finsler seulement $C^0$, ce qui ne permet pas sans 
autre d'adapter l'in\'egalit\'e de Cheeger riemannienne, d'autre 
part, il y a diff\'erentes fa\c{c}ons d'induire une mesure sur les 
hypersurfaces selon que l'on s'int\'eresse \`a l'in\'egalit\'e 
isop\'erim\'etrique ou \`a la formule de la co-aire (comme souvent ces 
notions co\"{i}ncident d\`es que l'on se trouve dans un contexte 
riemannien). Aussi avons nous trait\'e ce probl\`eme d'un point du vue 
g\'en\'eral des vari\'et\'es de Finsler aux sections 3 (pour la 
formule de la co-aire) et 4 (pour l'in\'egalit\'e de Cheeger). Ces 
deux sections sont int\'eressantes pour elles-m\^emes, mais peuvent 
\'egalement \^etre admises pour la preuve du th\'eor\`eme \ref{equiv1}

\section{La non nullit\'e du bas du spectre implique la $\delta$-hyperbolicit\'e}

\begin{theo}\label{lambdaIdelta}
Un convexe $(\mathcal{C},d_\mathcal{C})$ du plan, muni de sa m\'etrique 
de Hilbert et
dont le bas du spectre est non nul est Gromov hyperbolique. Plus 
pr\'ecis\'ement :
$$
\forall\lambda>0, \exists\delta>0 \text{ tel que } \forall\mathcal{C}, \lambda_1(\mathcal{C}) \geq \lambda \Rightarrow (\mathcal{C},d_\mathcal{C}) \text{ est } %
\delta\text{-hyperbolique.}  
$$
\end{theo}

Soit $G_n:={\rm PGL}(\R^{n+1})$, $\pP^n:=\pP(\R^{n+1})$ l'espace projectif de $\R^{n+1}$.

Une partie \textsl{proprement convexe} $\mathcal{C}$ de $\pP^n$ est une partie 
convexe dont l'adh\'erence est incluse
dans le compl\'ementaire d'un hyperplan projectif. Ainsi on peut lui associer
un convexe born\'e de $\R^n$ et, inversement, \`a un convexe born\'e de $\R^n$ on peut associer
une partie proprement convexe de $\pP^n$.
Dans la suite, on 
d\'esignera par
$X_n$ l'ensemble des ouverts proprement convexes 
et par $X_n^\delta$ l'ensemble des ouverts proprement 
convexes et
$\delta$-hyperboliques de $\mathbb{P}^n$, munis de leur m\'etrique de Hilbert.

Nous munissons $X_n$ de la distance de Hausdorff entre les ensembles, 
comme d\'efinie dans \cite{benoist}, p. 2.

On peut alors \'enoncer:

\begin{prop}[\cite{benoist} proposition 2.11] \label{propbenoist}
Soit $F$ un sous ensemble ferm\'e de $X_n$, $G_n$-invariant, dont tous les \'el\'ements sont
strictement convexes (c'est \`a dire que l'int\'erieur du segment reliant deux 
points du bord de $\mathcal C$ est dans $\mathcal C$). Alors il existe 
un nombre r\'eel $\delta>0$ tel que $F\subset X_n^\delta$.
\end{prop}

En sorte que le th\'eor\`eme \ref{lambdaIdelta} 
se d\'eduit ais\'ement de la
proposition \ref{propbenoist} gr\^ace \`a la proposition suivante.
\begin{prop}\label{propferme}
Soit $\lambda>0$ et $F_\lambda\subset X_n$ l'ensemble des convexes de $X_n$ dont
le bas du spectre 
est sup\'erieur ou \'egal \`a $\lambda$. Alors $F_\lambda$ est une partie 
ferm\'ee et
$G_n$-invariante de $X_{n}$, et si $n=2$, tous
ses \'el\'ements sont strictement convexes.
\end{prop}

L'id\'ee de la preuve de cette proposition est la suivante: 

\begin{itemize}

\item La $G_n$-invariance est imm\'ediate, car pour tout
$g\in G_n$ et $\mathcal{C}\in X_n$, $\mathcal{C}$ et $g\mathcal{C}$ sont 
des espaces m\'etriques isom\'etriques, et ils ont donc le m\^eme bas 
du spectre.

\item La difficult\'e principale r\'esidera dans la fermeture.
Supposons qu'il existe une suite $(\mathcal{C}_i)$ dans $F_\lambda$ 
telle que $\mathcal{C}_i$ converge vers $\mathcal{C}$ au sens de 
Hausdorff
dans $X_n$ et $\lambda_1(\mathcal{C})< \lambda$. Alors, pour tout 
$\varepsilon >0$, on va montrer 
qu'il existe un domaine $\Omega$ tel que 
$\overline{\Omega} \subset \mathcal C$ avec $\lambda_{1}^{\mathcal 
C}(\Omega) \leq \lambda_{1}(\mathcal C) + \varepsilon$. Par ailleurs, pour 
$i$ assez grand, on aura $\overline{\Omega} \subset \mathcal C_{i}$. On en 
d\'eduira que $\lambda_{1}(\mathcal C_{i}) < \lambda$, ce qui est une 
contradiction. 

\item Le fait que la limite $\mathcal C$ est strictement convexe 
d\'ecoulera alors directement. 

\end{itemize}

On formalise cette id\'ee intuitive \`a l'aide des deux 
lemmes suivants.

 \begin{lemme}\label{firstlema} Soit $\mathcal C\subset\mathbb R^n$ un domaine convexe. 
Alors, pour tout ${\varepsilon>0}$, il existe un domaine $\Omega$ non vide avec
$\overline{\Omega} \subset \mathcal C$ tel que 
$$
\lambda_{1}^{\mathcal C}(\Omega) \leq \lambda_{1}(\mathcal C)+\varepsilon\text{.}
$$
\end{lemme}

\begin{proof} Fixons un point origine $P_{0} \in 
\mathcal C$.

Par d\'efinition de $\lambda_{1}(\mathcal 
 C)$, il existe une fonction $f$ \`a support dans $\mathcal C$ de 
quotient de Rayleigh inf\'erieur \`a $\lambda_{1}(\mathcal C) + \varepsilon$. 

   La distance euclidienne entre le support de $f$ et le bord 
$\partial \mathcal C$ de $\mathcal C$ \'etant strictement positif, 
il existe $\alpha <1$ tel que
 $f$ est encore \`a support dans un homoth\'etique $\alpha \mathcal 
 C$ de $\mathcal C$ par rapport au point $P_{0}$.

On choisit alors $\Omega = \alpha \mathcal C$, et, toujours par 
d\'efinition, le bas du spectre de $\Omega$ (comme domaine de 
$\mathcal C$) est major\'e par le quotient de Rayleigh de $f$, donc 
par $\lambda_{1}(\mathcal C) + \varepsilon$.
\end{proof}  

\noindent
\textbf{Notation:} Soit $\mathcal C$ un domaine convexe et $P_{0}\in 
\mathcal C$. Pour $\rho >0$, on va consid\'erer les convexes 
$(1+\rho)\mathcal C$ et $(1-\rho)\mathcal C$ qui sont les 
 homoth\'etiques de $\mathcal C$ par une homoth\'etie de centre 
$P_{0}$ et de rapport $(1+\rho)$ et $(1-\rho)$ respectivement.

\begin{lemme}\label{secondlema} Soit $\mathcal C$ un domaine convexe
      et $\Omega$ un 
domaine tel que $\overline{\Omega} \subset \mathcal C$. Soit $\varepsilon 
>0$. Alors, il existe $0< \rho <1$ tel que $\overline{\Omega} \subset 
(1-\rho)\mathcal C$ et tel que pour tout convexe $\Gamma$ 
v\'erifiant $(1-\rho)\mathcal C \subset \Gamma \subset (1+\rho) 
\mathcal C$, on ait
$$(1-\varepsilon) \lambda_{1}^{\mathcal C}(\Omega) \leq
\lambda_{1}^{\Gamma}(\Omega) \leq (1+\varepsilon) \lambda_{1}^{\mathcal 
C}(\Omega).$$
\end{lemme}

\begin{proof} Comme au lemme pr\'ec\'edent, on montre 
qu'il existe $\rho_{0} >0$ avec 
$$ \overline{\Omega} \subset (1-2 \rho_{0}) 
\mathcal C \subset (1-\rho_{0})\mathcal C \text{.}$$ 
Ainsi, pour tout 
$0<\rho<\rho_{0}$, la distance euclidienne entre $\overline{\Omega}$ et 
$\partial\bigl((1-\rho)\mathcal C\bigr)$ est uniform\'ement minor\'ee par une 
 constante strictement positive.

 D\`es maintenant, on suppose $\rho < \rho_{0}$ et $(1-\rho)\mathcal{C}\subset\Gamma\subset(1+\rho)\mathcal{C}$.

  Soit $p \in \Omega$ et $u_{p} \in T_{p}\mathcal C$, un vecteur non nul. Comme cons\'equence directe 
de la d\'efinition de la norme de Finsler associ\'ee \`a la 
distance de Hilbert, on a
\begin{equation}
\Vert u_{p} \Vert_{(1+\rho)\mathcal C} \leq \Vert u_{p} \Vert_{\Gamma} 
\leq \Vert u_{p} \Vert_{(1-\rho)\mathcal C}\text{.}
\end{equation}
On en d\'eduit d'une part l'in\'egalit\'e suivante, par dualit\'e, pour une forme lin\'eaire non nulle $l_{p}$ 
sur $T_{p}\mathcal C$:
 $$\Vert l_{p} \Vert_{(1-\rho)\mathcal C}^{*} 
\leq \Vert l_{p} \Vert_{\Gamma} ^{*}
\leq \Vert l_{p} \Vert_{(1+\rho)\mathcal C}^{*},$$
d'autre part la relation d'inclusion v\'erifi\'ee par les boules 
unit\'es de l'espace tangent en $p$:
    $$TB_{(1-\rho)\mathcal C}(p) \subset TB_{\Gamma}(p) \subset 
TB_{(1+\rho)\mathcal C}(p).$$
Ce qui implique que les densit\'es de volume associ\'ees 
v\'erifient
\begin{equation}
\frac{\omega_{n}}{\vole\bigl(B_{\mathcal (1+\rho) C}(p)\bigr)} \leq 
\frac{\omega_{n}}{\vole\bigl(B_{\Gamma}(p)\bigr)} \leq \frac{\omega_{n}}{\vole\bigl(B_{(1-\rho)\mathcal C}(p)\bigr)}.      
\end{equation}

De plus, par compacit\'e de $\Omega$ et par continuit\'e, les rapports 
\begin{align}
\frac{\Vert u_{p} \Vert_{(1+\rho)\mathcal C}}
{\Vert u_{p} \Vert_{(1-\rho)\mathcal C}}&, & \frac{\Vert l_{p} \Vert_{(1-\rho)\mathcal C}^{*}}
{\Vert l_{p} \Vert_{(1+\rho)\mathcal C}^{*}} && \text{et }&&\frac{\vole(B_{\mathcal (1+\rho) C}(p))}{\vole(B_{(1-\rho)\mathcal C}(p))}
\end{align} 
sont uniform\'ement contr\^ol\'es en fonction de $\rho$ par des 
fonctions tendant vers $1$ lorsque $\rho \to 0$.
Cela \'etant ind\'ependant du choix de $\Gamma$, on peut comparer 
deux convexes quelconques du type de $\Gamma$, soit ici $\Gamma$ 
et $\mathcal C$, et on obtient
 l'existence d'une
 fonction $h=h(\rho)$ telle que $h(\rho) \to 0$ lorsque $\rho 
\to 0$ et telle que pour une fonction $f$ de 
classe $C^{\infty}$ \`a support dans $\Omega$, on ait 
$$\bigl(1-h(\rho)\bigr) \lambda_{1}^{\mathcal C}(\Omega) \leq
\lambda_{1}^{\Gamma}(\Omega) \leq \bigl(1+h(\rho)\bigr) \lambda_{1}^{\mathcal C}(\Omega).
$$
Il suffit donc de choisir $\rho$ assez proche de $0$ pour conclure.
\end{proof}

\begin{proof}[D\'emonstration de la proposition \ref{propferme}.]  
Pour voir que $F_{\lambda}$ est ferm\'e, on con\-sid\`ere une suite $(\mathcal C_{i})_{i\in {\mathbb N}}$ 
dans $F_{\lambda}$ qui converge vers un convexe $\mathcal C$ 
au sens de Hausdorff. Supposons que $\mathcal C \not \in 
F_{\lambda}$, c'est-\`a-dire qu'il existe $\varepsilon >0$ avec 
$\lambda_{1}(\mathcal C) \leq \lambda-\varepsilon$. Alors, par le 
premier lemme, il existe un domaine $\Omega$ avec $\overline{\Omega} 
\subset \mathcal C$ et 
$$\lambda_{1}^{\mathcal C}(\Omega) \leq  (\lambda - \varepsilon)+\varepsilon/3= \lambda-2\varepsilon/3.$$ 

Le deuxi\`eme 
lemme implique l'existence de $\rho >0$ tel que pour tout convexe $\Gamma$ 
avec $(1-\rho)\mathcal C \subset \Gamma \subset (1+\rho)\mathcal 
C $, on ait 
$$ \lambda_{1}^{\Gamma}(\Omega) \leq (1+\varepsilon/2)\lambda_{1}^{\mathcal C}(\Omega) \leq (1-\varepsilon/6) \lambda < \lambda \text{.}$$ 

Cela est en particulier vrai pour les 
\'el\'ements $\mathcal C_{i}$ de la suite lorsque $i$ est assez grand 
et contredit le fait que $\mathcal C_{i} \in F_{\lambda}$.

Il reste alors \`a montrer que tous les convexes de $F_\lambda$ 
sont strictement convexes en dimension deux. Cependant
si un convexe n'est pas strictement convexe, il admet un triangle dans
l'adh\'erence de son orbite (voir \cite{benoist} corollaire 2.9,
d\^u \`a Benz\'ecri \cite{benzecri}), et on obtient une contradiction comme dans la 
premi\`ere partie, car le bas du spectre d'un triangle
est nul. On montre ce dernier point par un calcul direct, \`a cause 
de la simplicit\'e de la g\'eom\'etrie de Hilbert du triangle. Si on 
consid\`ere une boule $B_{R}$ de rayon $R$, on voit que son aire est de 
l'ordre de $R^2$ et la longueur de son bord de l'ordre de $R$ lorsque 
$R \to \infty$. On construit alors une fonction test valant $1$ sur 
$B_{R}$, $0$ hors de $B_{R+1}$, et valant $1-d(p,B_{R})$ pour $p \in 
B_{R+1}-B_{R}.$ Le quotient de Rayleigh d'une telle fonction est de 
l'ordre de $1/R$ lorsque $R \to \infty$.
\end{proof}

\begin{rem}
 En ce qui concerne les triangles, une autre mani\`ere de voir que le bas du spectre
est nul, est d'utiliser le th\'eor\`eme 5 dans \cite{vernicos}, qui g\'en\'eralise l'in\'egalit\'e
de Faber-Krahn aux espaces vectoriels norm\'es de dimension finie. Celui-ci permet de dire que le
$\lambda_1$ d'une boule de rayon $R$ est de l'ordre de $1/R^2$, puisque la g\'eom\'etrie
du triangle est isom\'etrique \`a celle du plan muni de la norme dont la boule unit\'e
est un hexagone r\'egulier (voir \cite{dlharpe}).  
\end{rem}

En fait, la proposition $\ref{propferme}$ implique le th\'eor\`eme suivant,
valable en toute dimension.

\begin{theo}
Soit $(\mathcal{C},d_\mathcal{C})$  un convexe de $\R^n$ muni de sa m\'etrique de Hilbert
tel qu'il existe un r\'eel strictement positif $\lambda$ minorant la premi\`ere valeur propre
de toute section planaire, i.e. pour tout $H$ plan affine de $\R^n$ tel que $\mathcal{C}\cap H\neq \emptyset$
on a $$ \lambda_1(\mathcal{C}\cap H)\geq \lambda \text{.}$$
Alors il existe $\delta(\lambda)>0$ tel que $(\mathcal{C},d_\mathcal{C})$ est $\delta$-hyperbolique.
\end{theo}

\begin{proof}
 L'ensemble des sections planaires de $\mathcal{C}$ est un sous-ensemble de l'espace $F_\lambda$ 
de la proposition \ref{propferme}. Ainsi toute section planaire de  $\mathcal{C}$ est 
$\delta$-hyperbolique pour un m\^eme $\delta$ suivant la proposition \ref{propbenoist} appliqu\'ee \`a $F_\lambda$. 
Tout triangle \'etant dans une section planaire, cela termine la d\'emonstration.
\end{proof}

\section{De la $\delta$-hyperbolicit\'e \`a la non nullit\'e du bas du 
spectre}
Le but de cette section est la d\'emonstration du th\'eor\`eme suivant :
\begin{theo}\label{deltaIvol1} Pour tout r\'eel $\delta>0$, il existe $\lambda>0$ tel que, pour tout
$(\mathcal{C},d_h)$ convexe de $\R^n$ muni de sa m\'etrique de Hilbert on a:
Si $(\mathcal{C},d_h)$ est $\delta$-hyperbolique, alors pour tout $H$, plan affine de
$\R^n$ tel que $H\cap\mathcal{C}\neq\emptyset$, on a 
$$\lambda_1(\mathcal{C}\cap H)\geq \lambda \text{.}$$
\end{theo}
En dimension deux cela se traduit par:
\begin{cor}\label{deltaIvol}
Un convexe du plan $\mathcal{C}$ qui, muni de sa m\'etrique de Hilbert $d_\mathcal{C}$ est Gromov-hyperbolique,
  a son bas du spectre strictement positif. Plus pr\'ecis\'ement : 
  pour tout $\delta >0$, il existe $\lambda>0$
tel que, pour tout convexe du plan $C$:
$$
(\mathcal{C},d_\mathcal{C}) \text{ est } \delta\text{-hyperbolique} \Rightarrow \lambda_1(\mathcal{C})\geq\lambda\text{.}
$$
\end{cor}
Un point essentiel de la preuve est de contr\^oler le volume des 
boules de $(\mathcal C,d_{\mathcal C})$, et nous traitons ce point 
pour lui-m\^eme.

\subsection{Volume des boules}\label{volboule}
On se propose ici de montrer qu'en g\'e\-o\-m\'e\-trie de Hilbert, le 
volume des boules est contr\^ol\'e par leur rayon, et cela 
ind\'ependemment de la g\'eom\'etrie du convexe consid\'er\'e.

 \noindent
\textbf{Notation:} Dans la suite, puisqu'il n'y aura plus 
d'ambigu\"\i t\'e sur le convexe que l'on consid\`ere, on indicera plus 
les boules par le convexe $\mathcal C$ dans lequel elles se trouvent, mais par 
leur rayon.

\begin{theo}\label{bruno} Il existe deux constantes positives $C_{1}(R,n),C_{2}(R,n)$ 
telles que si $\mathcal C$ est un ouvert convexe born\'e de $\mathbb 
R^n$ muni de sa distance de Hilbert usuelle $d_{\mathcal C}$ et de sa 
 mesure finsl\'erienne associ\'ee $\mu_{\mathcal C}$, alors, pour toute boule $B_{R}$ de 
 rayon $R$ dans $(\mathcal C, d_{\mathcal C})$, on a
$$\underbrace{\frac{\omega_n}{4e^{2nR}} \Bigl(\frac{e^{2R}-1}{e^{2(R+1)}-1}\Bigr)^n}_{C_1(R,n)}%
\leq \mu_\mathcal{C}(B_{R}) \leq \underbrace{\left(\frac{e^{4R}-1}{2}\right)^n \omega_n}_{C_2(R,n)}.$$
\end{theo}
On constate en particulier que lorsque $R \to 0$, le volume d'une 
boule de rayon $R$ dans $\mathcal C$ est de l'ordre de $R^n$, \`a 
des constantes multiplicatives pr\`es ne d\'ependant que de la 
dimension. En particulier, le volume d'une boule de rayon fix\'e 
est contr\^ol\'e ind\'ependamment de la position de son centre et 
de $\mathcal C$.

La premi\`ere \'etape de la preuve consiste \`a comparer le volume 
euclidien d'une boule de rayon $R$ dans $\mathcal C$ centr\'ee en $x$ 
et que l'on note $B_{R}(x)$ et le volume euclidien de la boule de 
rayon $R$ tangente au point $x$
(pour la m\'etrique d Finsler $F_{\mathcal C}$ de classe  $C^0$) qui sera 
not\'ee $TB_{R}(x)$.

\begin{prop}\label{pbruno} Il existe deux constantes positives 
$c_{1}(R,n)$, $c_{2}(R,n)$ 
(ind\'ependantes de $\mathcal C$ et de $x$) telles 
que 
$$\underbrace{\left(\frac{e^{2R}-1}{2Re^{2R}}\right)^n}_{c_1(R,n)}%
\leq \frac{\vole B_{R}(x)}{\vole TB_{R}(x)} 
\leq \underbrace{\left(\frac{e^{2R}-1}{2R}\right)^n}_{c_2(R,n)}$$
o\`u $\vole$ d\'esigne le volume euclidien.
\end{prop}

\begin{proof}[D\'emonstration de la proposition \ref{pbruno}.] On regarde
 se qui se passe dans 
une direction donn\'ee, le long d'une droite. Pour all\'eger les 
notations et simplifier les calculs, on va travailler en dimension 1 
et supprimer autant que possible les normes intervenant dans 
l'expression de la m\'etrique de Finsler.

\begin{figure}[htpb]
\includegraphics{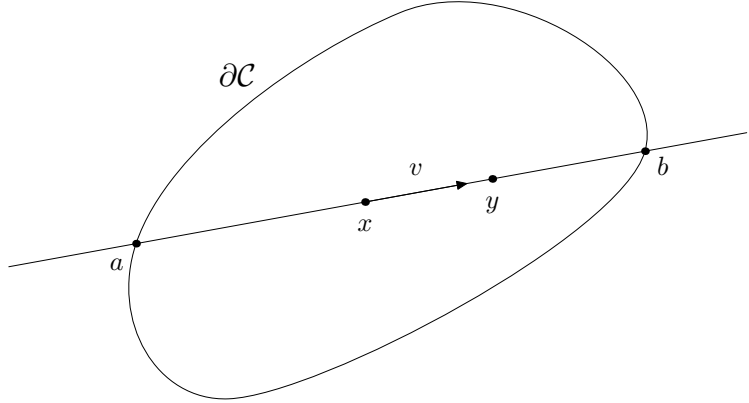}
\caption{Comparaison de $\Vert v \Vert_{e}$ et $y-x$}
\end{figure}

On se trouve sur une droite avec trois points 
ordonn\'es $a,x,b$. Un point $y>x$ tel que $d_{\mathcal C}(x,y)=R$ 
est tel que

 \begin{equation}
\label{eqdc}
\frac{1}{2} \ln (\frac{y-a}{x-a}\frac{b-x}{b-y}) =R
\end{equation}
et un vecteur $v \in T_{x}\mathcal C$ pointant vers $b$ et de norme 
de Finsler ${\Vert v \Vert_{\mathcal C}=R}$ sera tel que
\begin{equation}
\label{eqnormf}
\frac{1}{2} \Vert v \Vert_{e} (\frac{1}{b-x} + \frac{1}{x-a}) =R. 
\end{equation}

On cherche \`a comparer $\Vert v \Vert_{e}$ et $y-x$.

Commen\c cons par \'etudier $y-x$. 

Puisque $d_{\mathcal C}(x,y)=R$, en prenant l'exponentielle de l'\'equation (\ref{eqdc}) on obtient
 $$\frac{y-a}{x-a}\frac{b-x}{b-y} =e^{2R},$$
puis en isolant $y$
$$y = \frac{a(b-x)+ e^{2R}b(x-a)}{(b-x)+e^{2R}(x-a)},$$
on en d\'eduit alors une expression de $y-x$ en fonction de $a$,$b$,$x$ et $R$:
 $$y-x = \frac{(x-a)(b-x)(e^{2R}-1)}{(b-x)+e^{2R}(x-a)}.$$
On exprime finalement $y-x$ en fonction  des nouvelles variables ${\alpha=x-a}$, ${\beta=b-x}$ et ${\gamma=b-a}$ :
\begin{equation}
\label{eqymx}
y-x = \frac{\alpha \beta (e^{2R}-1)}{\beta + e^{2R}\alpha} = 
\frac{\alpha \beta (e^{2R}-1)}{\gamma + (e^{2R}-1)\alpha}.
\end{equation}

En ce qui concerne $\n{v}_e$, en combinant l'expression de la norme de Finsler (\ref{eqnormf})
et le fait que $\Vert v \Vert_{\mathcal C}=R$, on obtient en fonction de $\alpha$,$\beta$,$\gamma$ et $R$:

\begin{equation}
\label{eqnormv}
\Vert v \Vert_{e} = \frac{2R(x-a)(b-x)}{b-a}= \frac{2R\alpha\beta}{\gamma}.
\end{equation}

 On utilise \`a pr\'esent les expressions (\ref{eqnormv}) de $\n{v}_e$ et (\ref{eqymx}) de $y-x$ 
pour obtenir une \'equation de leur rapport en fonction de $\alpha$,$\beta$,$\gamma$ et $R$: 
\begin{equation}
\label{eqrapport}
\frac{\Vert v \Vert_{e} }{y-x}= \frac{2R\alpha 
\beta}{\gamma}\frac{\gamma + (e^{2R}-1)\alpha}{\alpha \beta 
(e^{2R}-1)}=
\frac{2R}{(e^{2R}-1)} \frac{\gamma + (e^{2R}-1)\alpha}{\gamma}.  
\end{equation}
On remarque alors qu'on a l'in\'egalit\'e 
$$1 \leq \frac{\gamma + (e^{2R}-1)\alpha}{\gamma} \leq e^{2R},$$
 et on l'applique \`a l'\'equation (\ref{eqrapport}) pour trouver 
$$\frac{2R}{(e^{2R}-1)} \leq \frac{\Vert v \Vert_{e} }{y-x} \leq \frac{2R 
e^{2R}}{(e^{2R}-1)}.$$

 Cette derni\`ere in\'egalit\'e induit les inclusions suivantes
$$\frac{2R}{(e^{2R}-1)} B_{R}(x) \subset TB_{R}(x) \subset  \frac{2R 
 e^{2R}}{(e^{2R}-1)} B_{R}(x).$$

 Les volumes euclidiens respectifs de $B_{R}(x)$ et $TB_{R}(x)$ 
sont donc contr\^ol\'es par les constantes annonc\'ees dans 
l'\'enonc\'e de la proposition, et qui ne d\'ependent que de $R$ 
et de $n$.
\end{proof}
\begin{cor}\label{cortreiz} Comme la boule $TB_{R}(x)$ est homoth\'etique de rapport 
$R$ avec la boule $TB_{1}(x)$, on d\'eduit l'existence de deux 
constantes positives $c_{1}'(R,n)$ et $c_{2}'(R,n)$ telles que
$$
\underbrace{\left(\frac{e^{2R}-1}{2e^{2R}}\right)^n}_{c_1'(R,n)}\leq \frac{\vole B_{R}(x)}{\vole TB_{1}(x)}%
\leq \underbrace{\left(\frac{e^{2R}-1}{2}\right)^n}_{c_2'(R,n)}  
.$$
\end{cor}

\begin{proof}[D\'emonstration du th\'eor\`eme \ref{bruno}.] On cherche \`a contr\^oler le volume 
de la boule $B_{R}(x)$. On sait qu'il est d\'efini par
\begin{equation}
 \label{eqvolboule}
\mu_{\mathcal C}\bigl(B_{R}(x)\bigr) = \int_{B_{R}(x)} \frac{w_{n}}{\vole TB_{1}(y)}dy
\end{equation}
o\`u $dy$ d\'esigne la mesure de Lebesgue et $\omega_{n}$ le volume 
euclidien de la boule unit\'e.

On sait par ailleurs que pour tout $y \in B_{R}(x)$,
on a $B_{R}(x) \subset B_{2R}(y)$. En combinant ceci avec le corollaire 
\ref{cortreiz}, on obtient d'abord  
$$\frac{1}{\vole TB_{1}(y)} \leq c_{2}'(2R,n)\frac{1}{ \vole B_{2R}(y)} \leq 
c_{2}'(2R,n) \frac{1}{\vole B_{R}(x)},$$
et puis, par int\'egration,
$$\mu_\mathcal{C}\bigl(B_{R}(x)\bigr) \leq \omega_{n} c_{2}'(n,R) = C_{2}(R,n).$$
On a \'egalement $\vole TB_{R+1}(x) = (R+1)^n \vole TB_{1}(x)$. 
Comme $B_{1}(y) \subset B_{R+1}(x)$, suivant le corollaire \ref{cortreiz},
on a d'une part
\begin{multline*}
\frac{1}{\vole TB_{1}(y)}\geq \left(\frac{e^{2}-1}{2e^{2}}\right)^n\frac{1}{\vole B_{1}(y)}
\geq\frac{1}{4}\frac{1}{\vole B_{R+1}(x)},
\end{multline*}
et d'autre part
\begin{multline*}
\frac{1}{\vole B_{R+1}(x)}\geq\Bigl(\frac{2}{e^{2(R+1)}-1}\Bigr)^n\frac{1}{\vole TB_{1}(x)}\geq\\
\Bigl(\frac{2}{e^{2(R+1)}-1}\Bigr)^n\left(\frac{e^{2R}-1}{2e^{2R}}\right)^n%
\frac{1}{\vole B_{R}(x)}.
\end{multline*}
On obtient alors par int\'egration:
$$\mu_\mathcal{C}\bigl(B_{R}(x)\bigr) \geq \frac{\omega_n}{4e^{2nR}}%
\Bigl(\frac{e^{2R}-1}{e^{2(R+1)}-1}\Bigr)^n=C_{1}(R,n).$$
\end{proof}

\subsection{Constante de Cheeger et bas du spectre.}\label{Cheeger} Dans la suite, on va 
\'etudier la constante de Cheeger des sections planaires d'un 
convexe $\mathcal C$ muni de sa m\'etrique de Hilbert. Si $\mathcal C$ est un domaine convexe de $\mathbb R^2$,
la constante de Cheeger de $(\mathcal C,d_{\mathcal C})$ est d\'efinie par
$$
I_\infty(\mathcal{C})=\inf_{\Omega\subset \mathcal{C}}\frac{\nu_{\mathcal{C}}(\partial\Omega)}{\mu_{\mathcal{C}}(\Omega)}
$$
o\`u $\Omega$ est un domaine dont l'adh\'erence (prise dans $\mathbb 
R^2$) v\'erifie $\overline{\Omega} \subset \mathcal 
C$ et dont le bord $\partial \Omega$ est une courbe rectifiable de 
longueur $\nu_{\mathcal C}(\partial \Omega)$. 

On va d'abord montrer que la constante de Cheeger des sections planaires 
d'un convexe $(\mathcal C,d_{\mathcal C})$ hyperbolique au sens de 
Gromov est strictement positive.
\begin{prop}\label{delcheeger1}
Soit  $(\mathcal{C},d_h)$ un convexe de $\R^n$ muni de sa m\'etrique de Hilbert
qui est $\delta$-hyperbolique. 
Alors il existe une constante $I(\delta)>0$, d\'ependant de $\delta$, 
telle que pour tout plan affine $H$ de $\R^n$ avec 
$\mathcal{C}\cap H\neq \emptyset$, on a $I_\infty(\mathcal{C}\cap H)\geq I(\delta)$.
\end{prop}

En dimension deux, on a

\begin{cor}\label{delcheeger}
 Soit  $(\mathcal{C},d_h)$ un convexe du plan muni de sa m\'etrique de Hilbert qui est 
 $\delta$-hyperbolique. Alors il existe une constante $I(\delta)>0$, d\'ependant de $\delta$, telle que 
$I_\infty(\mathcal{C}) \geq I(\delta)$.
\end{cor}

Pour montrer cela on utilise un r\'esultat d\'ecrit par M.~Bridson et A.~Haefliger (cf. Chapitre III.H section 2,
d\'efinition 2.1 et proposition  2.7 de \cite{brihae}):

\begin{prop}[\cite{brihae}]
Soit $X$ un espace g\'eod\'esique. Si $X$ est $\delta$-hyperbolique, alors il existe deux constantes
$A$ et $B$ et $\varepsilon>0$ telles que pour toute courbe rectifiable 
$c: S^1 \to X$ on ait
$$
Aire_\varepsilon(c)\leq A\nu_{\mathcal{C}}(\partial\Omega) + B
$$
o\`u $Aire_\varepsilon(c)=\min\{|\Phi|\mid \Phi \text{ est un } \varepsilon-\text{ remplissage de }c\}$. 
\end{prop}

\begin{proof}[D\'emonstration de la proposition \ref{delcheeger1}] Sans limiter la 
g\'en\'eralit\'e, on se ram\`ene \`a la dimension 2. Pour un domaine 
$\Omega$ de $\mathcal C$ de bord rectifiable, on cherche \`a minorer 
le quotient $\nu_{\mathcal C}(\partial\Omega)/\mu_{\mathcal C}(\Omega)$. 
On peut d'une part supposer que $\Omega$ est connexe, mais \'egalement 
qu'il est simplement connexe. En effet, comme $\Omega$ est born\'e 
dans $\mathcal C$, son compl\'ementaire dans $\mathcal C$ compte 
exactement une composante connexe non born\'ee, du fait que $\mathcal 
C$ est hom\'eomorphe \`a $\mathbb R^2$. En 
consid\'erant le compl\'ementaire de cette composante connexe, on 
obtient un nouveau domaine born\'e et simplement connexe $\Omega'$ contenant $\Omega$, donc 
v\'erifiant $\mu_{\mathcal C}(\Omega') \geq \mu_{\mathcal C}(\Omega)$, 
et dont le bord $\partial \Omega'$ est contenu dans le bord de 
$\Omega$ et v\'erifie donc $\nu_{\mathcal C}(\partial \Omega') \leq \nu_{\mathcal 
C}(\partial \Omega)$. Ainsi, on a
$$\frac{\nu_{\mathcal C}(\partial \Omega')}{\mu_{\mathcal C}(\Omega')} \leq
\frac{\nu_{\mathcal C}(\partial \Omega)}{\mu_{\mathcal C}(\Omega)}.$$

Il suffit donc de faire la preuve pour les domaines simplement 
connexes.

Soit donc $\Omega$ un domaine simplement connexe dont le bord 
$\partial \Omega$ est une courbe rectifiable $c: S^1 \to \mathcal C$. 

Alors, puisque nous sommes en
dimension deux, chaque partie du $\varepsilon$-remplissage peut \^etre recouvert par une boule 
de rayon $\varepsilon$.  Pour pr\'eciser cela, rappelons qu'un 
$\varepsilon$-remplissage est la donn\'ee d'une triangulation $P$ du 
disque $D^2$ et d'une application (non n\'ecessairement continue) $\Phi$ 
de $D^2$  dans $\Omega$, qui applique le bord $S^1$ de $D^2$ sur le 
bord $\partial \Omega$ de $\Omega$, et telle que l'image de chaque 
face ferm\'ee de la triangulation $P$ a un diam\`etre au plus 
$\varepsilon$. Comme ici $\Omega$ est un disque topologique de $\mathbb 
R^2$, l'image des sommets de la triangulation doit \^etre 
$\varepsilon$-dense dans $\Omega$. En effet, on peut \'etendre 
affinement la restriction de $\Phi$ aux sommets de la 
triangulation en une application continue $\tilde{\Phi}$ du disque dans $\Omega$ dont 
la restriction au bord est \'egal \`a l'application $c$. Supposons que l'image 
des sommets de la triangulation ne soit pas $\varepsilon$-dense. Alors,
il existe 
un point $Q$ de $\Omega$ qui n'est pas dans l'image de 
$\tilde{\Phi}$. Cela provient du fait que l'image par $\tilde{\Phi}$ d'une face est le 
triangle d\'etermin\'e par l'image des sommets de la face 
consid\'er\'ee, qui est, par construction, de diam\`etre au plus 
$\varepsilon$. L'application $\tilde{\Phi}$ ne serait alors pas 
surjective, ce qui permettrait de r\'etracter le disque contin\^ument 
sur son bord.

Cela permet de recouvrir $\Omega$ 
par un nombre de boules de rayon $\varepsilon$ \'egal au nombre de faces 
de $P$.  

On en d\'eduit donc que
$$
\mu_{\mathcal{C}}(\Omega)\leq Aire_\varepsilon(c)C_2(\varepsilon,n)
\leq A'\nu_{\mathcal{C}}(\partial\Omega)+B' 
$$
et ainsi
$$\frac{\nu_{\mathcal{C}}(\partial\Omega)}{\mu_{\mathcal{C}}(\Omega)}\geq
\frac{1}{A'}-\frac{B'}{\mu_{\mathcal C}(\Omega)}.$$

Ceci nous donne une minoration de $I_\infty(\Omega)$ pour les domaines $\Omega$ de
volume $\mu_{\mathcal C}(\Omega) \geq 2A'B'$.
Pour les domaines tels $\mu_{\mathcal C}(\Omega) \leq 2A'B'$, on remarque que leur longueur
est au moins \'egale \`a
deux fois leur diam\`etre et que le volume d'une boule de rayon le diam\`etre,
contenant le domaine est plus grand que celui du domaine. Le th\'eor\`eme
\ref{bruno} permet de conclure que $I_\infty(\mathcal{C}\cap H)>0$. Remarquons que la minoration
de $I_\infty(\mathcal{C}\cap H)$ ne d\'epend que de $A$ et $B$. Ainsi dans un convexe de $\R^n$ dont la
g\'eom\'etrie de Hilbert est $\delta$-hyperbolique, la constante de Cheeger 
de toute section planaire est minor\'ee par une m\^eme constante.
\end{proof}

Il reste donc \`a en d\'eduire des informations pour le bas du spectre
$\lambda_1(\mathcal{C}\cap H)$. Classiquement, dans le cas des 
vari\'et\'es riemanniennes, cela se fait via l'in\'egalit\'e de 
Cheeger qui relie le bas du spectre et la constante de Cheeger. 
Cependant, dans le contexte des espaces m\'etriques, il n'est pas 
clair qu'une telle relation existe en toute g\'en\'eralit\'e. Aussi allons nous nous ramener 
dans un cadre o\`u on sait qu'une telle in\'egalit\'e existe: celui 
des g\'eom\'etries de Finsler. On sait que la distance de Hilbert 
d\'erive d'une m\'etrique de Finsler, mais cette derni\`ere a la 
r\'egularit\'e du bord du convexe consid\'er\'e $\mathcal C$. Le but 
est de d\'emontrer le th\'eor\`eme \ref{theofinal} ci-dessous, en se 
ramenant \`a une in\'egalit\'e de Cheeger classique dans un cadre 
finsl\'erien qui sera 
\'enonc\'ee au th\'eor\`eme \ref{chencheegertheo} et d\'emontr\'ee 
dans la section \ref{icheeg}.

\begin{theo}\label{theofinal}
Soit  $(\mathcal{C},d_{\mathcal C})$ un convexe du plan muni de sa m\'etrique de Hilbert. 
Alors il existe une constante $C>0$, ne d\'ependant pas de $\mathcal 
C$, telle que
$$\lambda_1(\mathcal{C})\geq C I^2_\infty(\mathcal{C})\text{.}$$
\end{theo}
On en d\'eduit le corollaire imm\'ediat suivant:
\begin{cor}\label{lambdacheeger}
Soit  $(\mathcal{C},d_{\mathcal C})$ un convexe du plan muni de sa
 m\'etrique de Hilbert. Si $\lambda_1(\mathcal C)=0$, alors $I_\infty(\mathcal C)=0$.
\end{cor}
Cela implique \'egalement le corollaire suivant, d\'ej\`a d\'emontr\'e
par A.~Karlsson et G.A.~Noskov \cite{kn} et 
par Y.~Benoist \cite{benoist} (voir Lemme 2.12 et corollaire 2.13).
\begin{cor}
Soit $(\mathcal{C},d_{\mathcal C})$ un convexe du plan muni de sa
m\'etrique de Hilbert qui est $\delta$-hyperbolique. Alors son bord
est $C^1$ et strictement convexe. 
\end{cor}

\begin{proof}
Dans le cas contraire, on utilise \`a nouveau le r\'esultat de 
  Benz\'ecri (voir \cite{benoist} corollaire 2.9,
d\^u \`a Benz\'ecri \cite{benzecri}), le convexe $\mathcal C$ admet un
triangle dans son adh\'erence, sous l'action de $G_2$. On a vu que cela 
impliquait 
le fait que $\lambda_1(\mathcal{C})=0$ (voir d\'emonstration de la proposition \ref{propferme}).

Le corollaire pr\'ec\'edent implique que la constante de Cheeger $I_\infty(C)$ est \'egalement nulle.
Ceci est en contradiction avec la proposition \ref{delcheeger1}.
\end{proof}

Remarquons enfin que l'on obtient \'egalement le th\'eor\`eme suivant
en dimension deux.

\begin{cor}
Soit $(\mathcal{C},d_{\mathcal C})$ un convexe du plan muni de sa
m\'etrique de Hilbert. Alors sa constante de Cheeger $I_\infty(C)$ est
non nulle si et seulement s'il est Gromov-hyperbolique.
\end{cor}
Ce corollaire est \`a comparer au r\'esultat de Cao \cite{jcao} dans le cas 
riemannien.

L'id\'ee de la preuve du th\'eor\`eme \ref{theofinal} est de remplacer l'\'etude de 
la relation entre le bas du spectre de $\mathcal C$, 
$\lambda_{1}(\mathcal C)$, et la constante de Cheeger de $\mathcal C$,
$I_{\infty}(\mathcal C)$, par l'\'etude des m\^emes objets pour des 
domaines ouverts $\Omega$ tels que $\overline{\Omega} \subset \mathcal 
C$. Cela permettra de se ramener \`a une in\'egalit\'e classique du 
type Cheeger dans le cas finsl\'erien, puis de conclure par passage 
\`a la limite.  

Le premier point est donc de d\'efinir la constante de Cheeger d'un 
domaine ouvert $\Omega$ de $\mathcal C$ avec $\overline{\Omega} \subset \mathcal 
C$: on pose
$$
I_\infty^{\mathcal{C}}(\Omega)=\inf_{U\subset \Omega}
\frac{\nu_{\mathcal{C}}(\partial U)}{\mu_{\mathcal{C}}(U)}
$$
o\`u $U$ est un domaine dont l'adh\'erence (prise dans $\mathbb R^2$) 
v\'erifie $\overline{U} \subset \Omega$
et dont le bord $\partial U$ est une courbe rectifiable de 
longueur $\nu_{\mathcal C}(\partial U)$. 

Le point d\'ecisif de la preuve \`a venir sera l'in\'egalit\'e de 
Cheeger qui est valable pour des m\'etriques de Finsler suffisamment 
r\'eguli\`eres. On utilise le r\'esultat suivant, qui est une application
dans notre cadre du th\'eor\`eme plus g\'en\'eral d\'emontr\'e dans la section
\ref{icheeg}.

\begin{theo}\label{chencheegertheo}
Soit  $(\mathcal{C},d_{\mathcal C})$ un convexe du plan muni de sa m\'etrique de Hilbert. 
Si le bord de $\mathcal{C}$ est $C^1$ et strictement convexe,
il existe une constante $C>0$ telle que pour tout ouvert $\Omega$ on a:
\begin{equation}
\label{chencheeger}
\lambda_{1}^{\mathcal{C}}(\Omega) \geq C\bigl(I_{\infty}^\mathcal{C}\bigr)^2(\Omega)\text{.}
\end{equation}
En particulier cela implique que
$$
\lambda_1(\mathcal{C})\geq C I^2_\infty(\mathcal{C})\text{.}
$$
\end{theo}
Remarquons que la constante $C$ est universelle; elle ne d\'epend 
en fait que de la dimension, comme on le verra au th\'eor\`eme 
\ref{n-chencheegertheo}.

\begin{proof}[D\'emonstration du th\'eor\`eme \ref{theofinal} \`a 
l'aide du th\'eor\`eme \ref{chencheegertheo}.] La seule difficult\'e, 
pour passer du th\'eor\`eme \ref{chencheegertheo} au th\'eor\`eme 
\ref{theofinal} est que dans ce dernier, le bord du convexe $\mathcal 
C$ n'est pas n\'ecessairement de classe $C^1$: il faut donc se 
ramener \`a cette situation.

Soit $\Omega$ un domaine du plan tel que ${\overline{\Omega}\subset \mathcal{C}}$. 
Comme dans les deux lemmes \ref{firstlema} et \ref{secondlema}, il existe des homoth\'etiques 
$(1+\rho) \mathcal C$ et $(1-\rho) \mathcal C$ de $\mathcal C$ tels 
que
$\overline{\Omega}\subset (1-2\rho)\mathcal C$, et tels que pour tout 
convexe $\Gamma$ v\'erifiant
$$(1-\rho) \mathcal C \subset \Gamma \subset(1+\rho) \mathcal C $$
on ait
$$\lim_{\rho \to 0}\lambda_{1}^{\Gamma}(\Omega)= 
\lambda_{1}^{\mathcal C}(\Omega).$$
Mais la preuve des lemmes \ref{firstlema} et \ref{secondlema} montre que la m\^eme chose est vraie pour 
la constante de Cheeger, c'est-\`a-dire que pour tout domaine 
$\Gamma$ v\'erifiant
$$(1-\rho) \mathcal C \subset \Gamma \subset(1+\rho) \mathcal C$$
on a
$$\lim_{\rho \to 0}I_{\infty}^{\Gamma}(\Omega) = 
I_{\infty}^{\mathcal C}(\Omega).$$

On fait alors tendre $\rho$ vers $0$ tout en choisissant une famille 
$\Gamma_{\rho}$ d'ensembles convexes de classe $C^1$, 
strictement convexes, tels que
$$(1-\rho) \mathcal C \subset \Gamma_{\rho} \subset(1+\rho) \mathcal C.$$

Pour montrer qu'une telle famille existe, on remplace dans un 
premier temps le convexe $\mathcal C$ par un polygone convexe 
$\mathcal C'$ approximant $\mathcal C$ de sorte \`a v\'erifier 
$$(1-\frac{\rho}{2}) \mathcal C \subset \mathcal C' 
\subset(1+\frac{\rho}{2}) \mathcal C,$$
puis on remplace chaque segment de $\mathcal C'$ par un arc 
strictement convexe, de fa\c{c}on \`a avoir un recollement de classe 
$C^1$ en chaque sommet. Cela peut se faire en approximant $\mathcal 
C'$ d'aussi pr\^et que l'on veut, \`a condition de faire 
l'approximation avec des fonctions polynomiales de degr\'e \'elev\'e.

D\`es lors, on peut appliquer l'in\'egalit\'e (\ref{chencheeger}) 
$$\lambda_{1}^{\Gamma_{\rho}}(\Omega) \geq C\bigl(I_{\infty}^{\Gamma_{\rho}}\bigr)^2(\Omega)$$
et, en passant \`a la limite, on obtient
$$
\lambda_{1}^{\mathcal C}(\Omega)\geq C\bigl(I_{\infty}^{\mathcal C}\bigr)^2(\Omega).
$$

Pour conclure, il suffit donc de passer \`a la limite sur les $\Omega$; 
si l'in\'egalit\'e ci-dessus n'\'etait pas vraie pour $\mathcal C$ \`a 
la place de $\Omega$, il existerait un domaine $\Omega$ pour lequel 
elle ne serait pas valide non plus, puisque la constante de Cheeger et 
le bas du spectre de $\mathcal C$ sont approximables d'aussi pr\^et 
que l'on veut par la constante de Cheeger et 
le bas du spectre des domaines 
d'adh\'erence compacte dans $\mathcal C$, voir lemme \ref{firstlema}.
\end{proof}

\begin{nb}
Dans notre preuve nous approchons par une famille de convexes dont le bord est de classe
$C^1$ et strictement convexe, ce qui n\'ecessite l'adaptation \`a ce cadre de l'in\'egalit\'e
(\ref{chencheeger}).
Un r\'esultat de L.~H\"ormander \cite{hormander}, pr\'ecis\'ement le lemme 2.3.2, permet
d'obtenir une famille d'ensemble de classe $C^2$, dont le bord admet un hessien strictement
positif en tout point. Dans ce cadre l'in\'egalit\'e (\ref{chencheeger}) est classique,
voir par exemple Z.~Shen \cite{zshen}.
\end{nb}

Le th\'eor\`eme \ref{deltaIvol1} est \`a pr\'esent un
corollaire du th\'eor\`eme \ref{theofinal} qui contr\^ole le bas du 
spectre $\lambda_{1}(\mathcal C)$ par la constante de Cheeger 
$I_{\infty}(\mathcal C)$
et de la proposition \ref{delcheeger1}
qui contr\^ole la constante de Cheeger des sections planaires en fonction du $\delta$ de la 
$\delta$-hyperbolicit\'e.

\section{Formule de la co-aire}
Afin de d\'emontrer une version tr\`es g\'en\'erale du th\'eor\`eme \ref{chencheegertheo},
on a besoin d'une formule de la co-aire.
Ce qui suit
apporte les outils et notions n\'ecessaires \`a sa d\'emonstration et son 
expression: dans une vari\'et\'e de Finsler $(M,F)$ munie d'une 
mesure $\mu$, il s'agit d'une part de 
d\'efinir une mesure naturellement associ\'ee \`a $\mu$ sur les 
hypersurfaces de $M$ (sous-section \ref{mesurehyp}), puis de 
pr\'eciser la notion de gradient (sous-section \ref{grad}). A la 
section \ref{icheeg}, on appliquera cela pour la mesure de Hausdorff 
associ\'ee \`a la m\'etrique de Finsler de la vari\'et\'e.

Cette section combine des id\'ees contenues dans les livres
de A.~Thompson \cite{thompson} chapitre 3 et chapitre 5 
et Z.~Shen \cite{zshen} chapitres 3 et 4.

\subsection{Mesure sur les hypersurfaces}\label{mesurehyp}
Soit $(V,F)$ un espace vectoriel de dimension finie muni d'une norme 
$F$, dont la sph\`ere unit\'e $S_F(0,1)$ centr\'ee en
$0$ est $C^1$ et strictement convexe.

\begin{defi}
On dira que $x$ est normal \`a $y$, ce que l'on notera $x \dashv y$, si et seulement si, pour
tout $\alpha\in \R$ on a
$$
F(x+\alpha y)\geq F(x)
$$
\end{defi}
\begin{rmq}
Cette relation n'est g\'en\'eralement pas sym\'etrique.
\end{rmq}

Soit $y$ un vecteur non nul. On lui associe l'hyperplan
$$
W_y := \left\{ w\in V \mid y\dashv w \right\}. 
$$

Consid\'erons une base $\{b_i\}_{i=2,\dots,n}$ de $W_y$ et $b_1=y$. Ainsi $\{b_i\}_{i=1}^n$ est une base de $V$.

D\'efinissons alors
\begin{gather*}
B^n(b_1,\dots,b_n):=\biggl\{ (y^i)_{i=1}^n\in \R^n, F\Bigr(\sum_{i=1}^ny^ib_i\Bigl)<1\biggr\} \quad \text{ et } \\
B^{n-1}_y(b_2,\dots,b_n):= \biggl\{ (y^j)_{j=2}^n\in \R^{n-1}, 
F\Bigr(\sum_{j=2}^ny^jb_j\Bigl)<1\biggr\}.
\end{gather*}

Les deux domaines $B^n$ et $B^{n-1}_y$ d\'ependent du choix de la base $\{b_j\}_{j=2}^n$.

On d\'efinit la fonction
$$
\zeta(y):=\frac{\vole({\mathbb B}^n)}{\vole(\mathbb B^{n-1})}
\cdot\frac{\vole\bigl(B^{n-1}_y(b_2,\dots,b_n)\bigr)}%
{F(y)\vole\bigl( B^n(b_1,\dots,b_n)\bigr)}
$$

Remarquons, que si $F$ est euclidienne, alors $\zeta(y)=1$.

De plus, dans tous les cas, puisque pour tout r\'eel positif $\lambda>0$, $F(\lambda y)=\lambda F(y)$,
\begin{gather*}
\vole\bigl(B^n(\lambda y,b_2,\dots,b_n)\bigr)=\vole\bigl(B^n(y,b_2/\lambda,\dots,b_n/\lambda)\bigr)\lambda^n \text{ et }\\ 
\vole\bigl(B^{n-1}(b_2,\dots,b_n)\bigr)=\vole\bigl(B^{n-1}(b_2/\lambda,\dots,b_n/\lambda)\bigr)\lambda^{n-1} 
\end{gather*}
on a la propri\'et\'e suivante:
$$
\zeta(\lambda y)=\zeta(y), \lambda>0, y\neq 0 \text{.}
$$
\begin{ppt}\label{pptequiv}
Soit $(V,F)$ un espace vectoriel norm\'e et 
$$
c_n:=\vole({\mathbb B}^n)/\vole({\mathbb B}^{n-1})\text{.}
$$
 Alors pour $y$ non nul on a
$$
2^{-n }c_n\leq\zeta(y)\leq2^n c_n\text{.}
$$
\end{ppt}

\begin{proof}
Soit $y=\sum_{i=1}^n y^ib_i$ un vecteur norm\'e de $(V,F)$ et consid\'erons la boule $B_y^{n-1}$.
Alors pour tout $t\in[0,1]$ et $(v^j)_{j=2}^n\in B_y^{n-1}$, par l'in\'egalit\'e triangulaire, on a
$$
F(ty+\sum_{j=2}^nv^jb_j)\leq tF(y)+F(\sum_{j=2}^nv^jb_j)\leq 2
$$
donc $[0,1]\times B_y^{n-1}\subset 2B^n$. De ceci on d\'eduit que 
$$
{\vole}_{n-1}(B_y^{n-1})\leq 2^n \vole(B^n)\text{.}
$$
Soit \`a pr\'esent $(t,(v^j)_{j=2}^n)\in B^n$. Autrement dit $F(ty+\sum_{j=2}^nv^jb_j)\leq 1$ 
et $\sum_{j=2}^nv^jb_j\in W_y$, soit $y\dashv \sum_{j=2}^nv^jb_j$. La normalit\'e donne ($t\neq 0$)
$$
|t|=F(ty)=F(|t|y)\leq F\bigl(|t|y+(|t|/t)\cdot \sum_{j=2}^nv^jb_j\bigr) \leq 1
$$
et, gr\^ace \`a l'in\'egalit\'e triangulaire, on obtient
$$
F(\sum_{j=2}^nv^jb_j)\leq F(ty+\sum_{j=2}^nv^jb_j)+F(-ty)\leq2
$$
ceci implique que $B^n\subset [-1,1]\times2[B_y^{n-1}]$, et donc
$$
\vole(B^n)\leq2^n {\vole}_{n-1}(B_y^{n-1})
$$
\end{proof}

On peut maintenant d\'efinir une mesure ad\'equate sur une 
hypersurface $N$ d'une vari\'et\'e de Finsler $(M,F)$. Sur $M$, on 
peut consid\'erer la mesure de Hausdorff associ\'ee \`a $F$, que 
l'on note $\mu_{F}$. Pour d\'efinir une mesure sur $N$, on 
consid\`ere la restriction de $F$ \`a $N$ qui muni $N$ d'une 
m\'etrique de Finsler dont la mesure de Hausdorff associ\'ee est 
not\'ee $\bar{\nu}_{F}$. Cette derni\`ere est en apparence le candidat 
naturel \`a \^etre la mesure associ\'ee \`a $\mu_{F}$ sur $N$. C'est 
d'ailleurs la mesure que l'on a consid\'er\'ee \`a la sous-section 
\ref{Cheeger} pour d\'efinir la constante de Cheeger 
$I_\infty(\mathcal{C})$ en dimension 2.
Cependant, la suite de ce paragraphe, notamment la preuve 
de la formule de la co-aire, montre qu'un choix un peu diff\'erent 
peut \^etre
plus ad\'equat. On introduit un champ de vecteurs normaux unit\'es $n$ sur 
$N$. Cela signifie que pour tout point $p \in N$, l'hyperplan $T_{p}N$ 
tangent \`a $N$ au point $p$ est normal \`a $n(p)$ au sens d\'efini 
ci-dessus, c'est-\`a-dire
$$
T_{p}N := \left\{ w\in T_{p}M \mid n(p)\dashv w \right\}.
$$
On choisit alors d'associer \`a $\mu_{F}$ la mesure 
$$\nu_{F} = \zeta (n) \bar{\nu}_{F}.$$

Ainsi, en g\'en\'eral, la mesure $\nu_{F}$ ne co\"{i}ncide pas avec 
$\bar{\nu}_{F}$, mais c'est le cas lorsque la m\'etrique est 
riemannienne. Remarquons que, par la propri\'et\'e \ref{pptequiv},
ces deux mesures $\bar{\nu}_{F}$ et $\nu_{F}$ sont dans un rapport 
born\'e, uniform\'ement contr\^ol\'e, ce qui sera crucial \`a la 
section \ref{icheeg}, puisque l'on utilisera $\bar{\nu}_{F}$ pour 
d\'efinir la constante de Cheeger et $\nu_{F}$ pour la formule de la 
co-aire.

Cette d\'efinition, pour le cas particulier de la mesure 
de Hausdorff, se g\'en\'eralise comme suit:
\begin{defi}\label{defmeshyp}
Soit $N$ une hypersurface dans un espace de Finsler $(M,F)$. Si $\mu$ 
est une mesure sur $M$ absolument continue par rapport \`a la mesure 
de Hausdorff $\mu_{F}$, i.e. on peut \'ecrire 
$\mu(p)=\phi(p)\mu_{F}(p)$ en tout point $p$ de $M$, on associe \`a 
$\mu$ la mesure $\nu$ sur l'hypersurface $N$ d\'efinie par
$$\nu(p)=\phi(p) \nu_{F}(p)$$
pour tout $p \in N$.
\end{defi}
Remarquons que l'on doit introduire la pond\'eration par la fonction 
$\zeta$ lorsqu'on cherche \`a obtenir la propri\'et\'e naturelle 
suivante.
\begin{ppt}
Si $V_F$ est le volume de Hausdorff induit par $F$ sur $M$, alors
$$  A_F\bigl(S(x,\rho)\bigr)=\lim_{\varepsilon\to 0}\frac{V_F\bigl(B(x,\rho+\varepsilon)\bigr)-V_F\bigl(B(x,\rho)\bigr)}{\varepsilon}$$
\end{ppt}

\subsection{Gradient d'une fonction par rapport \`a une norme strictement 
convexe}\label{grad}
A la diff\'erentielle $df_x$ d'une fonction en $x\in M^n$, 
on voudrait associer un vecteur $X$ dans l'espace tangent $T_xM$
tel que
$$
F^*(df_x)=F(X)
$$
Pour ceci, on prend le point $X\in S_F\bigl(0,F^*(df_x)\bigr)$, sur la sph\`ere de centre $0$ et de
rayon $F^*(df_x)$, tel que 
$$
\ker df_x =\bigl\{v \mid X \dashv v\bigr\} 
$$
et $df_x(X)>0$ (l'hyperplan $X+\ker df_x$ est l'hyperplan d'appui
de $S\bigl(0,F^*(df_x)\bigr)$ en X).

On remarquera qu'en raison de la stricte convexit\'e
$$
df_x(X)=F^*(df_x)F(X)=F^*(df_x)^2\text{.}
$$

On dira que $X$ est le gradient de $f$ en $x$ et on le notera $\nabla f(x)$.

Remarquons que pour une fonction $C^1$ la fonction $x\mapsto \nabla f(x)$ est continue et
que pour une fonction au moins $C^2$ elle est $C^1$, puisque la sph\`ere $S(0,1)$ est $C^1$.

\begin{lemme}[Lemme de Gauss Finsl\'erien]Soit $f$ une fonction $C^2$ sur un ouvert $U$ tel que $df\neq 0$. Soit $N_t=f^{-1}(t)$. Alors le
vecteur ${n=\nabla f_{N_t}/F(\nabla f)}$ est normal \`a $N_t$. Autrement dit le gradient est
normal aux lignes de niveaux.  
\end{lemme}

\begin{proof}
Puisque $f$ est constante sur la ligne de niveau $N_t$, on a que
$$
df(w)=0,\ \forall w\in TN_t
$$
et donc par d\'efinition de $\nabla f$, $\nabla f\dashv w,\ \forall w \in TN_t$. 
\end{proof}

\subsection{D\'emonstration de  la formule de la Co-aire}
Dans cette section on d\'emontre le th\'eor\`eme suivant.
\begin{theo}
Soit $(M^n,F,\mu)$ une vari\'et\'e finsl\'erienne mesur\'ee, telle qu'en tout point $x\in M$
la boule unit\'e de la norme $F$ est $C^1$ et strictement convexe. 
Soit $f$ une fonction $C^1$ par morceaux
sur $M$ tel que pour tout $t$, l'ensemble $f^{-1}(t)$ est compact. Alors pour toute fonction continue 
$\varphi$ sur $M$ on a la formule de la co-aire suivante
\begin{equation}
\label{coarea}
\int_M\varphi F^*(df)d\mu = \int_{-\infty}^{+\infty}\int_{f^{-1}(t)} \varphi 
d\nu\text{,}
\end{equation}
o\`u $\nu$ d\'esigne la mesure associ\'ee \`a $\mu$ sur 
l'hypersurface $f^{-1}(t)$ selon la d\'efinition \ref{defmeshyp}.
\end{theo}

On va faire une d\'emonstration partielle de ce r\'esultat classique, en g\'eom\'etrie riemannienne,
dans le cas o\`u les lignes de niveau 
$\{f^{-1}(t)\}$ sont des sous-vari\'et\'es, cela pour 
montrer comment la mesure $\nu$ que l'on a introduit sur les 
hypersurfaces entre naturellement en jeu. Par le th\'eor\`eme de Sard, 
on sait qu'il existe un ensemble de mesure nulle de points $t$ pour 
lesquels $\{f^{-1}(t)\}$ n'est pas une sous-vari\'et\'e.  On s'inspire de la d\'emonstration de \cite{zshen}. 

\begin{proof}
On va faire la preuve dans une carte $U$. Pour simplifier, on supposera que $df\neq 0$ sur $U$
et que $f$ est $C^2$.
Fixons un nombre $t_0$ tel que $f^{-1}(t_0)\cap U\neq 0$ et soit $X$ le champ de vecteurs d\'efini sur
$U$ par
$$
X=\frac{\nabla f}{F^*(df)^2}\text{.}
$$

On commence par relever le flot.
Pour un point $x\in f^{-1}(t_0)\cap U$, soit $c(t)$ la courbe int\'egrale de $X$ telle que $c(t_0)=x$
(on peut parler de \guillemotleft\,la\,\guillemotright{} courbe car $X$ est $C^1$ sur $U$).
Alors
$$
\frac{d}{dt}\bigl[ f\circ c(t)\bigr]=df(X)=\frac{df(\nabla f)}{F^*(df)^2}=1
$$
par construction de $\nabla f$ et par cons\'equent
$$
f\circ c(t) =t \text{.}
$$

Les courbes int\'egrales de $X$ donnent naissance \`a un syst\`eme de coordonn\'ees
$\psi=(x^1,\dots,x^n)\colon U\to\mathopen]-\varepsilon,\varepsilon\mathclose[\times{\mathbb B}^{n-1}$ tel que
$$
f\circ\psi^{-1}(x^1,(x^i)_{i=2}^n)=x^1\text{.}
$$
D\'efinissons l'ensemble $N_t$ par

$$
f^{-1}(t)\cap U=\psi^{-1}\bigl( {t}\times{\mathbb B}^{n-1}\bigr) \text{.}
$$

Suivant le lemme de Gauss finsl\'erien, le champ de vecteurs

$$
n=\frac{\nabla f}{F^*(df)}
$$
est normal \`a $N_t$. Consid\'erons donc une base locale $\{b_i\}_{i=1,\dots,n}$ de $TM$ avec
$b_1=n$ et $b_i=\frac{\partial}{\partial x^i}$ pour $i=2,\dots,n$. Soit $(\theta^i)_{i=1,\dots,n}$ la base
duale sur $T^*M$. On a alors :
\begin{align*}
 \theta^1=\lambda dx^1=\lambda df, && \theta^i=dx^i,\ (i=2,\dots,n)\text{.}
\end{align*}
Remarquons \'egalement que
$$
1=\theta^1(n)=\lambda df\biggl(\frac{\nabla f}{F^*(df)}\biggr)=\lambda F^*(df)\text{,}
$$
ce qui d\'etermine $\lambda$, \textsl{i.e.}
$$
\theta^1=\frac{1}{F^*(df)}dx^1 \text{.}
$$

A pr\'esent, si $d\mu=\sigma(x)\theta^1\land\dots\land\theta^n$, alors par construction
$$
d\nu=\sigma(x){\hat \theta}^2\land \dots \land {\hat \theta}^n,
$$
o\`u les ${\hat \theta}^i$ sont les tir\'es en arri\`ere des $\theta^i$ sur $N_t$.
Autrement dit
\begin{eqnarray}
d\mu&=&\frac{\sigma(x)}{F^*(df)}dx^1\land \dots \land dx^n ;\\
d\nu&=&\sigma(x) dx^2\land \dots \land dx^n.
\end{eqnarray}
\begin{nb}
Si $d\mu$ est la mesure de finsler $n$-dimensionnelle, alors
on voit ici appara\^\i tre naturellement la fonction $\zeta$ puisque
$d\nu$ est exactement $\zeta$ fois la mesure de finsler $n-1$ dimensionnelle.
\end{nb}

Il reste \`a int\'egrer sur $U$:

\begin{eqnarray*}  
\int_U\varphi F^*(df)d\mu &=&\int_{\mathopen]-\varepsilon,\varepsilon \mathclose[\times{\mathbb B}^{n-1}} f\sigma(x)\,dx^1\dots dx^n; \\
&=& \int_{-\varepsilon}^{\varepsilon}\left( \int_{{\mathbb B}^{n-1}} f\sigma(x) dx^2\dots dx^n\right) dx^1;\\
&=& \int_{-\varepsilon}^\varepsilon \left( \int_{N_t\cap U} fd\nu\right)dt \text{.}
\end{eqnarray*}

\end{proof}

\section{In\'egalit\'e de Cheeger Finsl\'erienne}\label{icheeg}

Dans cette section, on va d\'emontrer une version n-dimensionelle du 
th\'eor\`eme \ref{chencheegertheo}:

\begin{theo}\label{n-chencheegertheo}
Soit  $(\mathcal{C},d_{\mathcal C})$ un domaine convexe de $\mathbb R^n$ muni de sa m\'etrique de Hilbert. 
Si le bord de $\mathcal{C}$ est $C^1$ et strictement convexe,
il existe une constante $C=C(n)>0$ telle que pour tout ouvert $\Omega$ on a:
\begin{equation}
\label{chencheegerbis}
\lambda_{1}^{\mathcal{C}}(\Omega) \geq C\bigl(I_{\infty}^\mathcal{C}\bigr)^2(\Omega)\text{.}
\end{equation}
En particulier cela implique que
$$
\lambda_1(\mathcal{C})\geq C I^2_\infty(\mathcal{C})\text{.}
$$
\end{theo}
Ce th\'eor\`eme sera une cons\'equence imm\'ediate de l'in\'egalit\'e 
de Cheeger en g\'eom\'etrie de Finsler. Si $(M,F)$ est une 
vari\'et\'e de Finsler de classe $C^1$ de dimension $n$, $\mu_{F}$ la mesure de 
Hausdorff associ\'ee \`a $F$, on d\'esigne par 
$\lambda_{1}(M,F)$ le bas du spectre de $(M,F)$ d\'efini par
$$\lambda_{1}(M,F)= \inf_f \frac{\int_M F^*(df)^2d\mu_{F}}{\int_\Omega 
|f|^2d\mu_{F}},$$
o\`u $f$ est une fonction lipschitzienne \`a support compact dans $M$. 
On d\'esigne $I_{\infty}(M,F)$ la constante de Cheeger de $(M,F)$ 
d\'efinie par
$$I_\infty(M,F)=\inf_{U}\frac{\bar{\nu}_F(\partial U)}{\mu_F(U)} 
$$
o\`u $U$ d\'ecrit les ouverts de $M$ d'adh\'erence compacte et dont le bord 
$\partial U$ est une sous-vari\'et\'e de dimension $n-1$. On rappelle 
que $\bar{\nu}_{F}$ d\'esigne la mesure de Hausdorff des 
hypersurfaces associ\'ee \`a la restriction de $F$. Alors, dans ces 
conditions
\begin{theo}[In\'egalit\'e de Cheeger finslerienne] \label{cheegfins} Il existe une constante 
$C=C(n)>0$ telle que
$$\lambda_{1}(M,F) \geq C I_\infty^2(M,F).$$
\end{theo}

De mani\`ere g\'en\'erale, soit $(M^n,F,\mu)$ une vari\'et\'e finsl\'erienne 
mesur\'ee, et $\nu$ la mesure induite sur les hypersurfaces au sens de 
la d\'efinition \ref{defmeshyp}.

Consid\'erons un domaine $\Omega$ de $M$ d'adh\'erence compacte.
On d\'efinit la constante de Sobolev $S_\infty(\Omega)$ de $\Omega$ 
par :
$$
S_\infty(\Omega)=\inf_f \frac{\int_\Omega F^*(df)d\mu}{\int_\Omega |f|d\mu}
$$
pour $f$ une fonction lipschitzienne \`a support compact dans $\Omega$. 
On d\'efinit de mani\`ere similaire $S_\infty(M)$ par
$$
S_\infty(M)=\inf_f \frac{\int_M F^*(df)d\mu}{\int_M |f|d\mu}
$$
pour $f$ une fonction lipschitzienne \`a support compact dans $M$. 

Aux mesures $\mu$ et $\nu$, on peut associer une constante de Cheeger 
sur $\Omega$ et sur $M$ d\'efinie par
\begin{gather*}
I_\infty^{M,\mu}(\Omega)=\inf_{U}\frac{\nu(\partial U)}{\mu(U)}\\
\left(\text{resp. } I_\infty^{\mu}(M) =   \inf_{U}\frac{\nu(\partial U)}{\mu(U)}\right)\text{.}
\end{gather*}
o\`u $U$ est un domaine tel que $\overline{U} \subset \Omega$ (resp. $\overline{U} \subset M$)
et dont le bord $\partial U$ est une sous-vari\'et\'e de dimension 
$n-1$.

Comme on l'a dit, lorsque $\mu = \mu_{F}$, la mesure de Hausdorff 
associ\'ee \`a $F$, on a $I_\infty^{\mu_{F}}(M) \not = 
I_{\infty}(M,F)$ puisque $\nu_{F} \not = \bar{\nu}_{F}$.

L'in\'egalit\'e de Cheeger est une cons\'equence de
\begin{theo}\label{egalite}
    
On a les \'egalit\'es $I_\infty^{M,\mu}(\Omega)=S_\infty(\Omega)$ et $I_\infty^{\mu}(M)=S_\infty(M)$.
\end{theo}

En effet, en consid\'erant $f=h^2$ et en appliquant l'in\'egalit\'e de 
Cauchy-Schwarz, on a
$$\bigl(I_\infty^{\mu_{F}}\bigr)^2(M) \leq 4\inf_{h} \dfrac{\int_M 
F^*(dh)^2d\mu_{F}}{\int_M |h|^2d\mu_{F}} = \lambda_{1}(M,F).$$

Il reste alors \`a comparer $I_\infty^{\mu_{F}}(M)$ \`a la constante 
de Cheeger $I_{\infty}(M,F)$ que l'on consid\`ere dans la premi\`ere 
partie. Cela vient de l'in\'egalit\'e (\ref{pptequiv}) qui nous permet d'obtenir l'in\'egalit\'e
 \begin{equation}
\label{cheegequivresp}
I_\infty(M)c_1(n)\leq I_\infty^{\mu_F}(M)\leq I_\infty(M)c_2(n)\text{.}
\end{equation}

\begin{proof}[D\'emonstration du th\'eor\`eme \ref{egalite}]
On suit 
la d\'emonstration du th\'eor\`eme 6.2 dans \cite{chavel} page 266--269.
 
On va d\'emontrer $I_\infty^{\mu}(M)=S_\infty(M)$.
La d\'emonstration consiste \`a d\'emontrer deux in\'egalit\'es.

\noindent\textbf{Premi\`ere \'etape: } On montre $S_\infty(M)\leq I_\infty^{\mu}(M)$. 
Soit donc $U$ un domaine ouvert de $M$
d'adh\'erence compacte et bord lisse. Pour tout $\varepsilon$ suffisamment petit on consid\`ere la fonction
$$f_\varepsilon(x)=
\begin{cases}
1& x\in U\\
1/\varepsilon d(x,\partial U) & x\in M\backslash U,\ d(x,\partial U)<\varepsilon\\
0 & x\in M\backslash U,\ d(x,\partial U)\geq\varepsilon
\end{cases}\text{.}
$$
en sorte que $f_\varepsilon$ est lipschitzienne. Donc
$$
S_\infty(M)\leq\frac{\int_M F^*(df_\varepsilon)d\mu}{\int_M|f_\varepsilon|d\mu}\text{.}
$$

Par construction, on obtient
$$
\lim_{\varepsilon\to0}\int_M|f_\varepsilon|d\mu = \mu(U)\text{.}
$$

De plus, on a
$$
F^*(df_\varepsilon) =
\begin{cases}
 1/\varepsilon & x\in M\backslash\overline{U},\ d(x,\partial U)\leq\varepsilon,\\
0 & sinon.
\end{cases}
$$
ce qui implique (voir section \ref{mesurehyp})
$$
\lim_{\varepsilon\to 0}\int_M F^*(df_\varepsilon)d\mu = \lim_{\varepsilon\to 0}\frac{\mu\bigl(\{ x\not\in U \:  d(x,\partial U)\leq \varepsilon\} \bigr)}{\varepsilon}=\nu(\partial U),
$$
et par cons\'equent, on en d\'eduit que pour tout $U$
$$
S_\nu(M)\leq\frac{\nu(\partial U)}{\mu(U)}
$$
ce qui permet de conclure cette premi\`ere \'etape en prenant 
l'infimum sur les $U$.

\noindent\textbf{Seconde \'etape:} On montre $I_\infty^{\mu}(M)\leq S_\infty(M)$.
Soit $f$ une fonction lisse \`a support compact dans $M$.
On pose 
\begin{align*}
U(t)=\bigl\{ x\: |f|(x)>t \bigr\},&& \mu(t)=\mu\bigl(U(t)\bigr), && \nu(t)=\nu\bigl(\partial U(t)\bigr).
\end{align*} 
pour les $t$ atteints par $f$. Alors la formule de la co-aire (\ref{coarea}) implique
$$
\int_MF^*(df)d\mu =\int_0^{+\infty} \nu(t)dt\geq I_\infty^{\mu}(M)\int_0^{+\infty} \mu(t) dt.
$$
(o\`u l'on a utilis\'e le fait que $F^*(df)=F^*(d|f|)$ presque partout sur $M$).

On remarque \'egalement que l'on a l'\'egalit\'e
$$\int_M |f|d\mu = \int_M \int_0^{|f|}dt d\mu =\int_0^{+\infty} \int_{U(t)}d\mu dt =\int_0^{+\infty}\mu(t)dt,
$$
qui permet d'arriver \`a l'in\'egalit\'e
$$
I_\infty^{\mu}(M)\leq \frac{\int_MF^*(df)d\mu}{\int_M |f|d\mu}
$$
et un passage \`a l'infimum sur toute les fonctions lisses \`a support compact $f$ permet de
conclure.
\end{proof}



































%% file: lamhyparxiv.bbl
\providecommand{\bysame}{\leavevmode\hbox to3em{\hrulefill}\thinspace}
\providecommand{\MR}{\relax\ifhmode\unskip\space\fi MR }
\providecommand{\MRhref}[2]{%
  \href{http://www.ams.org/mathscinet-getitem?mr=#1}{#2}
}
\providecommand{\href}[2]{#2}
\begin{thebibliography}{Tho96}

\bibitem[BBI01]{bbi}
D.~Burago, Y.~Burago, and S.~Ivanov, \emph{A course in metric geometry},
  Graduate Studies in Mathematics, vol.~33, American Mathematical Society,
  Providence, RI, 2001. \MR{2002e:53053}

\bibitem[Ben60]{benzecri}
J.P. Benz{\'e}cri, \emph{Sur les vari{\'e}t{\'e}s localement affines et
  localement projectives}, Bulletin de la Soci{\'e}t{\'e} Math{\'e}matique de
  France \textbf{88} (1960), 229--232.

\bibitem[Ben03]{benoist}
Y.~Benoist, \emph{Convexes hyperboliques et fonctions quasisym\'etriques},
  Publ. Math. Inst. Hautes \'Etudes Sci. (2003), no.~97, 181--237. \MR{2 010
  741}

\bibitem[BH99]{brihae}
M.R. Bridson and A.~Haefliger, \emph{Metric spaces of non-positive curvature},
  Comprehensive Studies in Mathematics, vol. 319, Springer, 1999.

\bibitem[Cao00]{jcao}
J.~Cao, \emph{Cheeger isoperimetric constants of {G}romov-hyperbolic spaces
  with quasi-poles}, Commun. Contemp. Math. \textbf{2} (2000), no.~4, 511--533.
  \MR{2001m:53068}

\bibitem[Cha93]{chavel}
I.~Chavel, \emph{Riemannian geometry: a modern introduction}, Cambridge Tracts
  in Mathematics, vol. 108, Cambridge University Press, Cambridge, 1993.
  \MR{95j:53001}

\bibitem[CV04]{cv}
B.~Colbois and P.~Verovic, \emph{Hilbert geometry for strictly convex domains},
  Geom. Dedicata \textbf{105} (2004), 29--42. \MR{2 057 242}

\bibitem[CVV]{cvv}
B.~Colbois, C.~Vernicos, and P.~Verovic, \emph{L'aire des triangles ideaux en
  g{\'e}om{\'e}trie de hilbert}, {{\`a}} para{{\^\i}}tre dans l'enseignement
  math{\'e}matique.

\bibitem[dlH93]{dlharpe}
P.~de~la Harpe, \emph{On {H}ilbert's metric for simplices}, Geometric group
  theory, Vol.\ 1 (Sussex, 1991), Cambridge Univ. Press, Cambridge, 1993,
  pp.~97--119.

\bibitem[H{\"o}r94]{hormander}
L.~H{\"o}rmander, \emph{Notions of convexity}, Progress in Mathematics, vol.
  127, Birkh\"auser Boston Inc., Boston, MA, 1994. \MR{MR1301332 (95k:00002)}

\bibitem[KN02]{kn}
A.~Karlsson and G.~A. Noskov, \emph{The {H}ilbert metric and {G}romov
  hyperbolicity}, Enseign. Math. (2) \textbf{48} (2002), no.~1-2, 73--89.
  \MR{2003f:53061}

\bibitem[She01]{zshen}
Zhongmin Shen, \emph{Lectures on finsler geometry}, World Scientific, 2001.

\bibitem[SM02]{so}
E.~Soci{\'e}-M{\'e}thou, \emph{Caract\'erisation des ellipso\"\i des par leurs
  groupes d'automorphismes}, Ann. Sci. \'Ecole Norm. Sup. (4) \textbf{35}
  (2002), no.~4, 537--548. \MR{1 981 171}

\bibitem[SM04]{so2}
\bysame, \emph{Behaviour of distance functions in {H}ilbert-{F}insler
  geometry}, Differential Geom. Appl. \textbf{20} (2004), no.~1, 1--10.
  \MR{2004i:53112}

\bibitem[Tho96]{thompson}
A.C. Thompson, \emph{Minkowski geometry}, Encyclopedia of Mathematics and its
  applications, vol.~63, Cambridge University Press, 1996.

\bibitem[Ver04]{vernicos}
C.~Vernicos, \emph{The macroscopic sound of tori}, Pacific J. Math.
  \textbf{213} (2004), no.~1, 121--156. \MR{MR2040254}

\end{thebibliography}
